\documentclass[11pt]{amsart}
\usepackage{geometry}                
\usepackage{graphicx}
\usepackage{amssymb}
\usepackage{epstopdf}
\usepackage{tikz}
\usetikzlibrary{arrows, decorations.markings, decorations.pathmorphing, backgrounds, positioning, fit, petri}
\usepackage{color}
\usepackage{float}

\newtheorem{thm}{Theorem}[section]
\newtheorem{prop}[thm]{Proposition} 
\newtheorem{cor}[thm]{Corollary} 
\newtheorem{lemma}[thm]{Lemma}

\theoremstyle{definition}
\newtheorem{example}{Example}
\newtheorem{defn}[thm]{Definition}
\newtheorem{remark}{Remark}

\title{Ribbon 2-knot groups of Coxeter type}
\author{Jens Harlander and Stephan Rosebrock}

\begin{document}

\maketitle

\centerline{Dedicated to the memory of Stephen Pride}

\begin{abstract}
Wirtinger presentations of deficiency 1 appear in the context of knots, long virtual knots, and ribbon 2-knots. They are encoded by (word) labeled oriented trees and, for that reason, are also called LOT presentations. These presentations are a well known and important testing ground for the validity (or failure) of Whitehead's asphericity conjecture. In this paper we define LOTs of Coxeter type and show that for every given $n$ there exists a (prime) LOT of Coxeter type with group of rank $n$. We also show that label separated Coxeter LOTs are aspherical.
\end{abstract}

\section{Introduction} 

Wirtinger presentations of deficiency 1 appear in the context of knots, long virtual knots, and ribbon 2-knots \cite{HarRose}. They are encoded by (word) labeled oriented trees and, for that reason, are also called LOT presentations.  Adding a generator to the set of relators in a Wirtinger presentation $P$ gives a balanced presentation of the trivial group. Thus the associated 2-complex $K(P)$ is a subcomplex of an aspherical (in fact contractible) 2-complex.  Wirtinger presentations are a well known and important testing ground for the validity (or failure) of Whitehead's asphericity conjecture, which states that a subcomplex of an aspherical 2-complex is aspherical. For more on the Whitehead conjecture see Bogley \cite{Bogley}, Berrick-Hillman \cite{BerrickHillman}, and Rosebrock \cite{Ro18}.

If $P$ is a Wirtinger presentation and the group $G(P)$ defined by $P$ is a 1-relator group, then $G(P)$ admits a 2-generator 1-relator presentation $P'$ and the 2-complex $K(P')$ is aspherical. Since $K(P')$ and $K(P)$ have the same Euler characteristic and the same fundamental group, it follows (using Schanuel's lemma and Kaplansky's theorem that states that finitely generated free $\mathbb ZG$-modules are Hopfian) that $K(P)$ is also aspherical. Thus, when investigating the asphericity of $K(P)$ for a given Wirtinger presentation $P$, the first thing to ask is if $G(P)$ is a 1-relator group.

Many knots have 2-generator 1-relator knot groups. Prime knots whose groups need more than 2 generators were known to Crowell and Fox in 1963. See Bleiler \cite{Bleiler} for a good discussion on this topic. As one example, Crowell and Fox consider a certain prime 9 crossing knot, show that its Wirtinger presentation simplifies to  
$$P=\langle x, y, z\ |\ y^{-1}xyx^{-1}y=x^{-1}zx^{-1}zxz^{-1}x, x^{-1}zxz^{-1}x=y^{-1}zyz^{-1}y \rangle,$$ and that the length of the chain of elementary ideals for this knot group is 2. It follows that the rank (=minimal number of generators) of $G(P)$ is greater than 2 and therefore equal to 3. This can also be seen without the use of elementary ideals. We have an epimorphism 
$$G(P)\to \Delta(3,3,3)=\langle x, y, z\ |\ x^2, y^2, z^2, (xy)^3, (xz)^3, (yz)^3 \rangle$$ sending $x\to x$,$y\to y$, $z\to z$. Since the rank of the Euclidian triangle group $\Delta(3,3,3)$ is 3 (see Klimento-Sakuma \cite{KlimentoSakuma}) we have $\mbox{rank}(G(P))=3$. 

This example is the motivation for this article. It is much easier to construct high rank ribbon 2-knot groups than classical knot groups, because we do not have to verify that a given Wirtinger presentation can be read off a knot projection (a 4-regular planar graph). Below we define (word) labeled oriented trees of Coxeter type and show that given a Coxeter group $W$, there exists a Coxeter type LOT group that maps onto $W$. Using this we give examples of LOT groups of arbitrarily high rank. 

In the second part of the paper we investigate the question of asphericity of LOTs of Coxeter type. We show that label separated LOTs of Coxeter type are aspherical. It turns out that the study of asphericity can be translated into questions concerning free subgroups of 1-relator LOT groups of dihedral type.

\vspace{0.5cm}

\section{Groups defined by graphs}
A {\em (word) labeled oriented graph} (LOG) is an oriented finite graph $\Gamma$ on vertices $\textbf{x}$ and edges ${\bf e}$, where each oriented edge is labeled by a word in ${\bf x}^{\pm 1}$. Associated with a LOG $\Gamma $ is the presentation 
$$P(\Gamma)=\langle {\bf x}\ |\  {\bf r}=\{r_e\ |\ e\in {\bf e} \} \rangle,$$ where 
$r_e=xw(wy)^{-1}$ in case $e=(x\stackrel{w}{\to}y)$ is the edge of $\Gamma$ starting at $x$, ending at $y$, and labeled with the word $w$ on letters in ${\bf x}^{\pm 1}$. We denote by $K(\Gamma)$ and $G(\Gamma)$ the standard 2-complex and the group defined by $P(\Gamma)$, respectively. The case where $\Gamma$ is a tree, now called a {\em (word) labeled oriented tree}  (LOT), is special. It is known that the groups $G(\Gamma)$, where $\Gamma$ is a LOT, are precisely the ribbon 2-knot groups (see Yajima \cite{Yajima}, Howie \cite{Howie}, and also Hillman \cite{Hillman}, section 1.7), since in that case $G(\Gamma)$ is a group of weight 1 (normally generated by a single element, in fact by each generator) that has a deficiency 1 presentation $P(\Gamma)$. The 2-complexes $K(\Gamma)$, $\Gamma$ a LOT,  are of central importance to Whitehead's asphericity conjecture, since adding a generator to the set of relators in $P(\Gamma)$ gives a balanced presentation of the trivial group. So $K(\Gamma)$ is a subcomplex of a 2-dimensional contractible complex. A question that has been open for a long time asks if $K(\Gamma)$ is aspherical, i.e.\ $\pi_2(K(\Gamma))=0$. See Bogley \cite{Bogley}, Berrick-Hillman \cite{BerrickHillman}, Rosebrock \cite{Ro18}.\\

Let $\Upsilon$ be a simplicial graph on vertices ${\bf x}$, and suppose edges $e$ are labeled with integers $m_e \ge 2$. Define
$$P(\Upsilon)=\langle {\bf x}\ |\ x^2, x\in {\bf x}, (xy)^{m_e}\ \mbox{if $e=\{ x,y\}$ is an edge} \rangle.$$ The group $W=W(\Upsilon)$ defined by this presentation is called a {\em Coxeter group}.
Let $K=K(\Upsilon)$ be the 2-complex associated with it. Consider the universal covering $\tilde K(\Upsilon)$. The 1-skeleton of $\tilde K(\Upsilon)$ is the Cayley graph for $(W,{\bf x})$. All edges in $\tilde K(\Upsilon)$ are double edges: For every $g\in W$ we have an edge $(g,x)$ connecting $g$ to $gx$, and an edge $(gx, x)$ connecting $gx$ to $g$. Note that a double edge pair bounds two 2-cells in $\tilde K(\Upsilon)$, coming from the relator $x^2$. We remove one and collapse the other one to an edge. This turns each double edge into a single unoriented edge.  Every relator $(xy)^{m_e}$ gives rise to $2m_e$ 2-cells with the same boundary. We remove all but one from this set. The 2-complex obtained in this fashion we denote by $\Sigma^{(2)}(\Upsilon)$. It is the 2-skeleton of the Coxeter complex $\Sigma(\Upsilon)$. See Proposition 7.3.4 in Davis \cite{Davis}. Under certain conditions, for example when ${\Upsilon}$ is a tree, the Coxeter complex is 2-dimensional: $\Sigma(\Upsilon)=\Sigma^{(2)}(\Upsilon)$. See \cite{Davis} Example 7.4.2. 

\begin{prop}\label{prop:treelike} Let $\Upsilon$ be a tree with associated Coxeter group $W(\Upsilon)$. Then
\begin{enumerate} 
\item For every edge $e=\{ x, y\}$ of $\Upsilon$ we have a 2-cell $\kappa_e$ in $\Sigma(\Upsilon)$ attached along a $2m_e$-gon whose edge labels read $(xy)^{m_e}$.
\item $\Sigma(\Upsilon)$ is the union of the 2-cells $w\kappa_e$, $e\in \{\mbox{edges of}\ \Upsilon\}$, $w\in W(\Upsilon)$. Furthermore, if $w_1\kappa_{e_1}\cap w_2\kappa_{e_2}\ne \emptyset$ then $e_1\cap e_2\ne \emptyset$; if $x=e_1\cap e_2$, then the edge $w_1\kappa_{e_1}\cap w_2\kappa_{e_2}$ carries the label $x$.
\item $\Sigma(\Upsilon)$ is a tree of 2-cells: If we connect the barycenters of the 2-cells with the barycenters of their boundary edges we obtain a tree. In particular, if $M$ is a finite connected union of Coxeter 2-cells $w_i\kappa_{e_i}$ in $\Sigma(\Upsilon)$, then there exists a 2-cell $w\kappa_e$ in $M$ that intersects with the rest of $M$ in a single edge.
\end{enumerate}
\end{prop}

\begin{figure}[htbp] 
   \centering
   \includegraphics[width=3in, height=5cm]{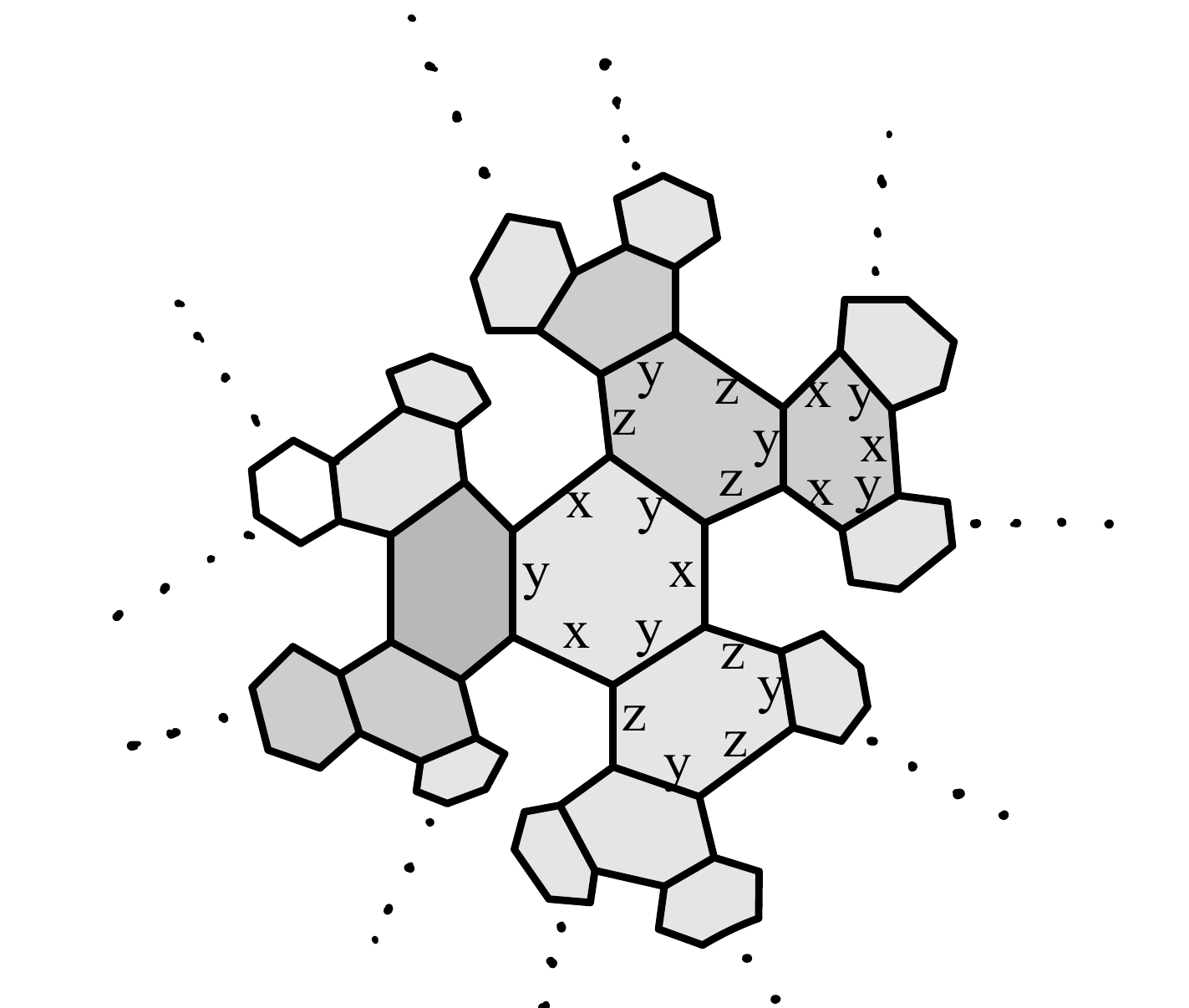} 
   \caption{The Coxeter complex $\Sigma(\Upsilon)$ for $\Upsilon=x\stackrel{3}{-}y\stackrel{3}{-}z$. It is a tree of Coxeter cells.}
   \label{fig:coxetertree}
\end{figure}

\noindent Proof. The statements (1) and (2) are clear from the construction of $\Upsilon$. 
For an edge $e=\{ x, y\}$ let $P(e)=\langle x, y\ | x^2, y^2, (xy)^{m_e} \rangle$. Let $D_{m_e}$ be the dihedral group defined by $P(e)$. Since $\Upsilon$ is a tree, $W(\Upsilon)$ is an amalgamated product of the $D_{m_e}$. The associated Bass-Serre tree can be seen inside the Coxeter complex $\Sigma(\Upsilon)$. The vertices of that tree are the barycenters of the 2-cells and 1-cells, and the edges connect barycenters of 2-cells to the barycenters of the 1-cells in the boundary of that 2-cell. We can think of $\Sigma(\Upsilon)$ as a tree of Coxeter 2-cells. An example is shown in Figure \ref{fig:coxetertree}. 

Suppose $L=\bigcup_{i=0}^{k}D_i$ is a union of 2-cells. Let $d_i$ be the barycenter of $D_i$. Let $d_p$ be a vertex in the Bass-Serre tree furthest away from $d_0$, $p\in \{ 0,\dots, k\}$. Consider a geodesic from $d_0$ to $d_p$ and let $d_q$ be the barycenter that is encountered just before getting to $d_p$ when traveling along the geodesic. Then $(\bigcup_{i\ne p} D_i)\cap D_p=D_q\cap D_p$, which is a single edge. \qed \\

The graph $\Upsilon$ also defines an {\em Artin presentation}. Denote by $prod(x,y,k)=xyxyx...$, where the length of the word is $k\ge 2$. Note that $prod(x,y,k)$ ends with $x$ in case $k$ is odd and it ends with $y$ if $k$ is even. Define
$$P_A(\Upsilon)=\langle {\bf x}\ |\ prod(x,y,m_e)=prod(y,x,m_e), \mbox{if $e=\{ x,y\}$ is an edge in $\Upsilon$} \rangle$$ 
and let $A(\Upsilon)$ be the group defined by $\Upsilon$. 

\begin{defn} Let $\Gamma$ be a LOT with vertex set ${\bf x}$. We say $\Gamma$ is of {\em Coxeter type} if for every edge $e=(x\stackrel{w}{\to}y)$ the word $w$ contains letters $z\ne x,y$ only with even (positive or negative) exponent.
\end{defn}



\begin{lemma}\label{lem:obvious} Let $\Gamma$ be a LOT of Coxeter type and $e=(x\stackrel{w}{\to}y)$ an edge. Then the relator $r_e=xw(wy)^{-1}$ reduces (up to cyclic permutation) to $\bar r_e=(yx)^{m_e}$, $m_e\ge 1$ and odd, in $\langle {\bf x}\ |\ x^2, x\in {\bf x} \rangle$.
\end{lemma}

\noindent Proof. The word $w$ reduced to an alternating word $\bar w$ in the letters $x$ and $y$. There are four cases to consider: 
\begin{enumerate}
\item $\bar w$ starts with $x$ and has even length; 
\item $\bar w$ starts with $x$ and has odd length; 
\item $\bar w$ starts with $y$ and has even length; 
\item $\bar w$ starts with $y$ and has odd length;
\end{enumerate}
In case (1) we have $\bar w=xyxy$, say. So $x(xyxy)y(xyxy)=xxyxyyxyxy=xy$. In case (2) we have $\bar w=xyxyx$, say. So $x(xyxyx)y(xyxyx)=xxyxyxyxyxyx=(yx)^5$. In case (3) we have $\bar w=yxyx$, say. So $x(yxyx)y(yxyx)=xy$. In case (4) we have $\bar w=yxyxy$, say. So $x(yxyxy)y(yxyxy)=(xy)^5$. \qed \\

Let $\Gamma$ be a LOT of Coxeter type. Define a tree $\Upsilon$ in the following way: Erase orientations in $\Gamma$ and if $e=(x\stackrel{w}{\to}y)$ is an edge and the LOT relator $r_e$ reduces to $\bar r_e=(yx)^{m_e}$ (up to cyclic permutation) in $\langle {\bf x}\ |\ x^2, x\in {\bf x} \rangle$, then label the (unoriented) edge $e$ by $m_e$.
We have a map $P(\Gamma)\to P(\Upsilon)$ sending $x$ to $x$ which induces a group epimorphism $G(\Gamma)\to W(\Upsilon)$. This process can be reversed. A LOT is {\em prime} if it does not contain a proper subLOT (proper means: not the entire LOT and not a single vertex).

\begin{lemma}\label{lem:GfromU} Let $\Upsilon$ be a Coxeter tree where all $m_e\ge 1$ and odd. Then there exists a (prime) LOT of Coxeter type $\Gamma$ such that the process just described produces $\Upsilon$ from $\Gamma$. In particular $G(\Gamma)$ maps onto $W(\Upsilon)$.
\end{lemma}

\noindent Proof. Suppose $e=\{ x, y\}$ is an edge in $\Upsilon$. Orient it in from $x$ to $y$. Let $w$ be a word in $\bf x$ that reduces to $(yx)^{\frac{m_e-1}{2}}$ in $\langle {\bf x}\ |\ x^2, x\in {\bf x} \rangle$. Let $e=(x\stackrel{w}{\to}y)$ be the corresponding edge in $\Gamma$. \qed 

\begin{remark} Let $\Upsilon$ be a Coxeter graph. If $e=\{ x, y\}$ is an edge  in $\Upsilon$, orient it from $x$ to $y$.  If $m_e$ is odd let $e=(x\stackrel{w_e}{\to}y)$ where $w_e=(yx)^{\frac{m_e-1}{2}}$. If $m_e$ is even remove the interior of the edge $e$ from $\Upsilon$ and attach a loop $e_x$ at the vertex $x$. Now $e_x=(x\stackrel{w_e}{\to}x)$ where $w_e=prod(y,x, m_e-1)$. Denote the labeled oriented graph we obtain by $\Gamma$. Observe that $P(\Gamma)=P_A(\Upsilon)$. Thus all Artin groups are LOG groups. This observation arose from a discussion with Gabriel Minian. 
\end{remark} 

\begin{example} Let $\Upsilon$ be a triangle with 3 vertices $x$, $y$, $z$ and edges $x\stackrel{3}{-} y$, $y\stackrel{3}{-} z$, $x\stackrel{2}{-} z$. We get a LOG $\Gamma$ with edges $x\stackrel{yx}{\to} y$, $y\stackrel{zy}{\to} z$, $x\stackrel{z}{\to} x$. Note that $G(\Gamma)=A(2,3,3)$ is an Artin group of spherical type with associated Coxeter group $W(2,3,3)$, the symmetric group $S_4$. It is known (Mulholland-Rolfsen \cite{MulhollandRolfsen}) that the commutator subgroup of $A(2,3,3)$ is finitely generated and perfect. So $G(\Gamma)=A(2,3,3)$ is not locally indicable. Whether LOT groups are locally indicable or not is (to our knowledge) an open problem.
\end{example}

\vspace{0.5cm}

\section{LOT groups of high rank}


\begin{thm}(Carette-Weidmann \cite{CaretteWeidmann})\label{thm:CW} Let $\Upsilon$ be a graph with $n$ vertices and assume that all the $m_{e} \ge 6\cdot 2^n$. Then the rank of $W(\Upsilon)$ is $n$.
\end{thm}

\begin{thm}\label{thm:Epi} Let $W=W(\Upsilon)$ be a Coxeter group such that $W_{ab}=\mathbb Z_2$. There exists a (prime) labeled oriented tree $\Gamma$ of Coxeter type so that $G=G(\Gamma)$ maps onto $W$. 
\end{thm}

\noindent Proof. Since $W_{ab}=\mathbb Z_2$ the Coxeter graph $\Upsilon$ is connected and contains a maximal tree $\Upsilon_0$ in which all $m_e$ are odd. Then $\Upsilon$ and $\Upsilon_0$ have the same set of vertices and we have an epimorphism $W(\Upsilon_0)\to W(\Upsilon)$. From Lemma \ref{lem:GfromU} we know that there is a (prime) LOT $\Gamma$ of Coxeter type so that $G(\Gamma)$ maps onto $W(\Upsilon_0)$. \qed

\begin{cor} For any given $n$ there exists a prime labeled oriented tree $\Gamma$ of Coxeter type with $n$ vertices so that $G(\Gamma)$ has rank $n$. In particular if $n\ge 3$ then $G(\Gamma)$ is not a 1-relator group. 
\end{cor}

\noindent Proof. This follows from Theorem \ref{thm:Epi} together with the Carette-Weidmann Theorem \ref{thm:CW}.  \qed

\begin{example}\label{ex:smallLOT} Let $\Gamma$ be the prime LOT $x\stackrel{yz^2x}{\rightarrow}y\stackrel{zx^2y}{\rightarrow}z$. Note that $G(\Gamma)$ maps onto the amalgamated product $D_3*_{\mathbb Z_2}D_3$ which can not be generated by two elements. Thus the rank of $G(\Gamma)$ is 3 and it follows that this LOT group is not a 1-relator group. If we drop the $z^2$ from the first edge word and $x^2$ from the second edge word we obtain a LOT $\Gamma_0$ that is not prime. In fact $G(\Gamma_0)=A(\Upsilon)$, where $\Upsilon$ is the Coxeter tree $x\stackrel{3}{-}y\stackrel{3}{-}z$.
\end{example}

\begin{remark} Note that if $\Gamma$ is a LOT of Coxeter type and $\Upsilon$ is the associated Coxeter tree, then $W(\Upsilon)$ is an amalgamated product of dihedral groups. A direct way to obtain upper bounds for the rank of $W(\Upsilon)$ without the full force of Theorem \ref{thm:CW} is via Weidmann \cite{Weidmann}. For example the LOT shown in Example \ref{ex:smallLOT} does not meet the conditions of the theorem.
\end{remark}

\begin{remark} A {\em reorientation} of a LOT is obtained when changing signs on the exponents of letters that occur in the edge words. Note that reorienting has no effect on the quotient $W(\Upsilon)$. Thus if $rk(G(\Gamma))=rk(W(\Upsilon))$, then this equation holds also for all reorientations of $\Gamma$.
\end{remark}

\vspace{0.5cm}

\section{Largeness}
 A group is {\em large} if it has a subgroup of finite index that has a free quotient of rank $\ge 2$. Large groups of deficiency 1 are studied in Button \cite{Button}. A list of properties can also be found in that paper. If $G$ is large then
 
 \begin{enumerate}
 \item $G$ contains free subgroups of rank $\ge 2$;
 \item $G$ is SQ-universal (every countable group is the subgroup of some quotient);
 \item $G$ has finite index subgroups with arbitrarily large first Betti number;
 \item $G$ has uniformly exponential word growth;
 \item $G$ has subgroup growth of strict type $n^n$ (which is the largest possible growth for finitely generated groups);
 \item the word problem for $G$ is solvable strongly generically in linear time.
 \end{enumerate}
 
 \begin{thm}\label{thm:large} Let $\Gamma$ be a LOT of Coxeter type on at least 3 vertices. Let $\Upsilon$ be the corresponding Coxeter tree and assume all $m_e\ge 3$. Then $G(\Gamma)$ is large.
 \end{thm}
 
 \noindent Proof. The conditions imply that $W(\Upsilon)$ is an infinite group that is a finite tree where the vertex groups are dihedral groups of type $m\ge 3$ and the edge groups are $\mathbb Z_2$. Thus $W(\Upsilon)$ contains a free subgroup $F$ of rank $\ge 2$ of finite index (see Serre's book Trees \cite{Serre}, Proposition 11, page 120). Let $H$ be the preimage of $F$ in $G(\Gamma)$. Then $H$ is a subgroup of $G(\Gamma)$ of finite index that maps onto $F$. It follows that $G(\Gamma)$ is large. \qed
 
 \begin{example} As in Example \ref{ex:smallLOT} let $\Gamma$ be the prime LOT $x\stackrel{yz^2x}{\rightarrow}y\stackrel{zx^2y}{\rightarrow}z$. Then $W(\Upsilon)=D_3*_{\mathbb Z_2}D_3$. Let $\Delta(3,3,2)$ be the spherical triangle group (it is the symmetric group $S_4$) defined by $Q=\langle  x, y, z \ |\ x^2, y^2, z^2, (xy)^3, (yz)^3, (xz)^2\rangle$. We have an epimorphism $W(\Upsilon)\to \Delta(3,3,2)$ and we claim that the kernel $V$ is free of rank $\ge 2$. Indeed, since both $D_3$'s of $W(\Upsilon)$ are also subgroups of  $\Delta(3,3,2)$, it follows that $V$ intersects both $D_3$'s trivially and it follows that $V$ acts freely of the Bass-Serre tree $T$ for $W(\Upsilon)=D_3*_{\mathbb Z_2}D_3$, and hence is free. Note that the valency of every vertex in $T$ is equal to 3 (since the index of $\mathbb Z_2$ in the $D_3$'s is 3), and so $V$ can not be cyclic. Here is why: Note that $V=\pi_1(X)$, where $X=T/V$ is a finite graph where every vertex has valency 3. So $v(X)=\frac{2e(X)}{3}$ and we obtain 
 $\chi(X)=v(X)-e(X)=\frac{2e(X)}{3}-e(X)<0$. Thus $\dim H_0(X)-\dim H_1(X)=1-\dim H_1(X)=\chi(X)<0$. So $\dim H_1(X)>1$ and hence $\dim V_{ab}>1$. One can also check directly that $(xz)^2$ and $x(xz)^2x^{-1}=(zx)^2$  generate a free subgroup of $V$ of rank 2.
 \end{example}
 
\newpage

\section{The question of asphericity}

\begin{figure}[htbp] 
   \centering
\begin{tikzpicture}[scale=0.94]

\draw (5,0) -- (5,-0.5) node[pos=0.5,right] {$z$};
\draw (6,0.5) -- (6.5,0.2)  node[pos=0.5,above,sloped] {$z$};
\draw (6,1.5) -- (6.5, 1.8)  node[pos=0.5,above,sloped] {$z$};
\draw (5,2) -- (5,2.5)  node[pos=0.5,right] {$z$};
\draw (4,1.5) -- (3.5,1.8)  node[pos=0.5,above,sloped] {$z$};
\draw (4,0.5) -- (3.5,0.2)  node[pos=0.5,above,sloped] {$z$};

\draw (5,0) -- (6, 0.5) node[pos=0.5,above] {$x$};
\draw (6, 0.5) -- (6,1.5)  node[pos=0.5,left] {$y$};
\draw (6,1.5) -- (5,2) node[pos=0.5,below] {$x$};
\draw (5,2) -- (4,1.5)  node[pos=0.5,below] {$y$};
\draw (4,1.5) -- (4,0.5)   node[pos=0.5,right] {$x$};
\draw (4,0.5) -- (5,0) node[pos=0.5,above] {$y$};

\draw[->] (2,1) -- (3,1);

\draw[blue] (0,0) -- (1, 0.5) node[pos=0.5,above] {$x$};
\draw[blue] (1, 0.5) -- (1,1.5)  node[pos=0.5,left] {$y$};
\draw[blue] (1,1.5) -- (0,2) node[pos=0.5,below] {$x$};
\draw (0,2) -- (-1,1.5)  node[pos=0.5,below] {$y$};
\draw (-1,1.5) -- (-1,0.5)   node[pos=0.5,right] {$x$};
\draw (-1,0.5) -- (0,0) node[pos=0.5,above] {$y$};

\draw[blue] (0,0) .. controls (0.5,0) .. (1, 0.5);
\draw[blue] (1, 0.5)  .. controls (1.3,1) .. (1,1.5);
\draw[blue] (1,1.5) .. controls (0.5,2) .. (0,2);
\draw (0,2)  .. controls (-0.5,2) ..  (-1,1.5);
\draw (-1,1.5)  .. controls (-1.2,1) .. (-1,0.5);
\draw (-1,0.5) .. controls (-0.5,0) ..  (0,0);

\draw (0,0) .. controls (0.2,-0.25) .. (0, -0.5)node[pos=0.5,right] {$z$};;
\draw (0,0) .. controls (-0.2,-0.25) .. (0, -0.5)node[pos=0.5,left] {$z$};;
\draw[blue] (1,0.5) .. controls (1.25,0.46) .. (1.5, 0.2)node[pos=0.5,right] {$z$};;
\draw[blue] (1,0.5) .. controls (1.25,0.2) .. (1.5, 0.2)node[pos=0.5,below] {$z$};;
\draw[blue] (1,1.5) .. controls (1.25,1.76) .. (1.5, 1.8)node[pos=0.5,above] {$z$};;
\draw[blue] (1,1.5) .. controls (1.25,1.56) .. (1.5, 1.8)node[pos=0.5,right] {$z$};;
\draw (0,2) .. controls (0.2,2.25) .. (0, 2.5)node[pos=0.5,right] {$z$};;
\draw (0,2) .. controls (-0.2,2.25) .. (0, 2.5)node[pos=0.5,left] {$z$};;
\draw (-1,1.5) .. controls (-1.2,1.7) .. (-1.5, 1.8)node[pos=0.5,above] {$z$};;
\draw (-1,1.5) .. controls (-1.3,1.5) .. (-1.5, 1.8)node[pos=0.5,left] {$z$};;
\draw (-1,0.5) .. controls (-1.3,0.5) .. (-1.5, 0.2)node[pos=0.5,above] {$z$};;
\draw (-1,0.5) .. controls (-1.3,0.2) .. (-1.5, 0.2)node[pos=0.5,right] {$z$};;

\end{tikzpicture}
   \caption{The complex $\bar K_e$ (on the left) in case $e=(x\stackrel{w}{\to}y)\in \Gamma$ with corresponding edge  $e=(x\stackrel{3}{\to}y)\in \Upsilon$, so the Coxeter relator is $(xy)^3$. On the right is the corresponding Coxeter cell $\kappa_e$ together with $z$-edges. The blue part is a $y$-side in $\bar K_e$.}
   \label{fig:Ke}
\end{figure}
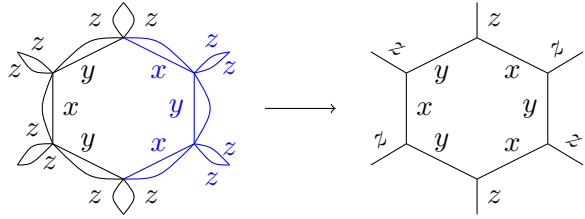

Let $\Gamma$ be a labeled oriented tree of Coxeter type and let $\Upsilon$ be the related Coxeter graph. Let $\bar K(\Gamma)$ be the normal covering space with fundamental group the kernel of the epimorphism $G(\Gamma)\to W(\Upsilon)$. We will analyze the structure of $\bar K(\Gamma)$. We have maps
$$\bar K(\Gamma)\to \tilde K(\Upsilon)\to \Sigma(\Upsilon),$$ and note that $\bar K(\Gamma)$ and $\tilde K(\Upsilon)$ have the same 1-skeleton. Let $e=(x\stackrel{w}{\to}y)$ be an edge in $\Gamma$. Let $P_e=\langle {\bf x}_e\ |\ r_e \rangle$, where ${\bf x}_e\subseteq {\bf x}$ is the subset of the vertices of $\Gamma$ that occur in $r_e$. Let ${\bf z}={\bf x}_e-\{ x, y\}$. Then $P_e=\langle x, y, {\bf z}\ |\ xw=wy \rangle$. The complex $K(P_e)$ is a subcomplex of $K(\Gamma)$. Consider the preimage of $K(P_e)$ under the covering projection $\bar K(\Gamma)\to K(\Gamma)$. It is a union of finite subcomplexes $w\bar K_e$, $w\in W(\Upsilon)$, that we will now describe in detail. The 1-skeleton of $\bar K_e$ is an $2m_e$-gon with double edges labeled in an alternating way by $x$ and $y$. At each of the $2m_e$ vertices we have a double edge for every $z\in {\bf z}$. The situation is depicted in Figure \ref{fig:Ke}. We have $2m_e$ 2-cells, attached along the loop with label $r_e$, starting at every vertex. 
The dihedral group $D_{m_e}$, the stabilizer of the Coxeter cell $\kappa_e$ in $\Sigma(\Upsilon)$, acts freely on $\bar K_e$. It is convenient to replace $\bar K_e$ with a complex with a single $D_{m_e}$ orbit of vertices. Let $\bar L_e$ be the 2 complex obtained from $\bar K_e$ in the following way: At every vertex collapse one of the $z$-edges from the $z$-double edge, $z\in {\bf z}$. The complex $\bar L_e$ is homotopy equivalent to $\bar K_e$. The 1-skeleton of $\bar L_e$ is an $2m_e$-gon with double edges labeled in an alternating way by $x$ and $y$. At each of the $2m_e$ vertices we have a loop for every $z\in {\bf z}$. Let $\hat r_e$ be the word obtained from $r_e$ by replacing every $z^p$, $z\in {\bf z}$, by $z^{\frac{p}{2}}$. Let $\hat P_e=\langle x, y, {\bf z}\ |\ \hat r_e \rangle$. Note that the dihedral group $D_{m_e}$ acts freely on $\bar L_e$ and we have a covering map $\bar L_e\to \bar L_e/D_{m_e}=K(\hat P_e)$.

\begin{lemma}\label{lem:aspherical} The 2-complex $\bar K_e$ is aspherical.
\end{lemma}

\noindent Proof. The complex $K(\hat P_e)$ is aspherical because $\hat P_e$ is a 1-relator presentation for which the relator is not a proper power. Thus $\bar L_e$ is aspherical, being a covering space of $K(\hat P_e)$. Since $\bar K_e$ is homotopy equivalent to $\bar L_e$, it follows that $\bar K_e$ is aspherical. \qed \\

An {\em $x$-side} of $\bar K_e$ consists of a double edge with label $x$ together with all the double edges connected to the two vertices of the $x$-double edge.  A $y$-side is defined in the same way. See Figure \ref{fig:Ke}, where the blue part on the left shows a $y$-side. Note that $\bar K_e$ has $m_e$ $x$-sides and $m_e$ $y$-sides. We refer to these as the {\em sides} of $\bar K_e$. We say $\bar K_e$ is {\em side injective} if the inclusion induced map $\pi_1(S)\to \pi_1(\bar K_e)$ is injective for every side $S$. An $x$-side in $\hat L_e$ is the image of an $x$-side under $\bar K_e\to \bar L_e$, etc.

\begin{lemma}\label{thm:sidebasis} Assume $m\ge 3$. Then $\bar K_e$ is $x$-side injective if and only 
$$\langle x^2, y^2, z, xy^2x^{-1}, xzx^{-1}, z\in {\bf z} \rangle$$ is a free subgroup of $G(\hat P_e)$ on the given basis.
\end{lemma}

\noindent Proof.  An $x$-side $S$ in $\bar L_e$ is an $x$-double edge, a $y$-double edge at each of the two vertices, and a loop for every $z\in {\bf z}$ at each of the two vertices. The image of $\pi_1(S)$ in $G(\hat P_e)$ under the covering projection is the group in the statement of the lemma. \qed \\

\begin{remark}Note that in case $m=1$ the situation can be different. The most extreme case occurs when $m_e=1$ for all edges $e$. In that case $W(\Upsilon)=\mathbb Z_2$ and so $\bar K(\Gamma)$ has only two vertices.  In this situation the 1-skeleton of each $\bar K_e$ is all of $\bar K(\Gamma)^{(1)}$, and in fact each side of $\bar K_e$ is all of $\bar K(\Gamma)^{(1)}$. So $\bar K_e$ is not side injective.
\end{remark}

\begin{lemma}\label{lem:subside} If $T$ is a subgraph of the 1-skeleton of $\bar K_e$ that does not involve every letter from ${\bf x_e}=\{ x, y, {\bf z}\}$, then $\pi_1(T)\to \pi_1(\bar K_e)$ is injective.
\end{lemma}

\noindent Proof. We can argue with $\bar L_e$ instead of $\bar K_e$. A reduced loop $\gamma$ in $T$ gives a reduced word $u$ in the generators of $\hat P_e$ that does not involve all letters from ${\bf x_e}=\{ x, y, {\bf z}\}$. The presentation $\hat P_e$ has only one relator $\hat r_e$ that does involve all letters from the generating set ${\bf x_e}=\{ x, y, {\bf z}\}$. The Freiheitssatz for 1-relator group implies that $u$ does not represent the trivial element in $G(\hat P)$. Thus $\gamma$ is not trivial in $\pi_1(\bar L_e)$. \qed \\

We continue our analysis. The complex $\bar K(\Gamma)$ is a union of the complexes $w\bar K_e$, $w\in W(\Upsilon)$, $e\in \mbox{edges of $\Gamma$}$. The maps
$$\bar K(\Gamma)\to \tilde K(\Upsilon)\to \Sigma(\Upsilon)$$
give a one-to-one correspondence between the $w\bar K_e$ and Coxeter cells $w\kappa_e$. Since $\Upsilon$ is a tree, the Coxeter complex $\Sigma(\Upsilon)$ is a tree of Coxeter cells $w\kappa_e$ and so $\bar K(\Gamma)$ is a tree of complexes $w\bar K_e$.  In complete analogy to Proposition \ref{prop:treelike} we have 

\begin{prop}\label{prop:alsotreelike} Consider $\bar K(\Gamma)=\bigcup w\bar K_e\to \Sigma(\Upsilon)=\bigcup w\kappa_e$.
\begin{enumerate} 
\item $\bar K(\Gamma)$ is the union of the 2-complexes $w\bar K_e$, $e\in \{\mbox{edges of}\ \Gamma\}$, $w\in W(\Upsilon)$. Furthermore, if $w_1\bar K_{e_1}\cap w_2\bar K_{e_2}\ne \emptyset$ then $e_1\cap e_2\ne \emptyset$; if $x=e_1\cap e_2$, then $w_1\bar K_{e_1}\cap w_2\bar K_{e_2}=T$, where $T$ is the subgraph of  an $x$-side $S$ that carries the letters ${\bf x}_{e_1}\cap {\bf x}_{e_2}$;
\item $\bar K(\Gamma)$ is a tree of 2-complexes. In particular, if $\bar M$ is a finite connected union of 2-complexes $w_i\bar K_{e_i}$ in $\bar K(\Gamma)$, then there exists a 2-complex $w\bar K_e$ in $\bar M$ that intersects with the rest of $\bar M$ in a subgraph of a single side.
\end{enumerate}
\end{prop}

\begin{thm}\label{thm:aspherical} Let $\Gamma$ be a LOT of Coxeter type. Then $K(\Gamma)$ is aspherical in case the $\bar K_e$ are side injective for every edge $e$ in $\Gamma$.
\end{thm}

\noindent Proof. We will show that $\bar K(\Gamma)$ is aspherical. It suffices to show that every finite union $\bar M=\bigcup_{i=1}^nw_i\bar K_{e_i}$ is aspherical.  We first claim that the sides of the $w_i\bar K_{e_i}$ $\pi_1$-inject into the union $\bar M$. We do induction on $n$. If $n=1$ the result follows from the hypothesis. Assume $n>1$. Then by Proposition \ref{prop:alsotreelike} part (2) there exists a 2-complex $w\bar K_e$ in $\bar M$ that intersects with the rest of $\bar M$ in a subgraph $T$ of a single side $S$ (of course $T$ could be $S$). Now by induction hypothesis the inclusion $S\subseteq \bar M-w\bar K_e=\bar M_0$ is $\pi_1$-injective and the inclusion $S\subseteq w\bar K_e$ is $\pi_1$-injective by hypothesis. It follows that $\pi_1(\bar M)$ is an amalgamated product $\pi_1(\bar M)=\pi_1(\bar M_0)*_{\pi_1(T)}\pi_1(w\bar K_e)$. And so the inclusion $S\subseteq \bar M$ is $\pi_1$-injective. All other sides that occur in $\bar M$ are either contained in $\bar M_0$ or in $w\bar K_e$. $\pi_1$-injectivity follows from the amalgamated product decomposition. Asphericity of $\bar M$ now follows from induction on $n$ and the amalgamated product decomposition $\pi_1(\bar M)=\pi_1(\bar M_0)*_{\pi_1(T)}\pi_1(w\bar K_e)$. \qed \\

\noindent {\bf Remark}. The above proof shows more than asphericity. Since each $\pi_1(\bar K_e)$ is a finite index subgroup of a 1-relator group, we see that $\pi_1(\bar K)$ is a tree of groups, the vertex groups being finite index subgroups of 1-relator groups, and the edge groups (over which we amalgamate) being finitely generated and free.

\begin{defn}A labeled oriented tree $\Gamma$ is called {\em label separated} if for every pair of edges 
$e_1$ and $e_2$ that have a vertex in common the intersection $\bf x_{e_1}\cap \bf x_{e_2}$ is a proper subset of both $\bf x_{e_1}$ and $\bf x_{e_2}$. \end{defn}

\begin{thm} Let $\Gamma$ be a label separated LOT of Coxeter type. Then $K(\Gamma)$ is aspherical.
\end{thm}

\noindent Proof. The proof is very much the same as the proof of Theorem \ref{thm:aspherical}. Let $\bar M=\bigcup_{i=1}^nw_i\bar K_{e_i}$as before. Again it suffices to show that $\bar M$ is aspherical. If $n=1$ that is clear. It is instructive to look at the case $n=2$. The intersection $w_1\bar K_{e_1}\cap w_2\bar K_{e_2}=T$ is the subgraph of a side that carries the letters ${\bf x}_{e_1}\cap {\bf x}_{e_2}$, which is a proper subset of both ${\bf x}_{e_1}$ and ${\bf x}_{e_2}$. $\pi_1$-injectivity for the inclusions $T\subseteq w_i\bar K_{e_i}$, $i=1,2$, follows from Lemma \ref{lem:subside}. We have $\pi_1(\bar M)=\pi_1(w_1\bar K_{e_1})*_{\pi_1(T)}\pi_1(w_2\bar K_{e_2})$ and $\bar M$ is aspherical. For $n\ge 2$ we argue by induction and obtain (as in the proof of Theorem \ref{thm:aspherical}) a decomposition $\pi_1(\bar M)=\pi_1(\bar M_0)*_{\pi_1(T)}\pi_1(w\bar K_e)$ which proves asphericity of $\bar M$. \qed

\vspace{0.5cm}

\section{Side injectivity}

Let $P=\langle a, b, {\bf c} \ |\ r \rangle$, be a 1-relator group, where $\bf c$ is a finite set of letters (which could be empty). We assume that $r$ is cyclically reduced and contains all generators. Assume further that $r=(ab)^m$ for some $m\ge 0$ modulo the relations $a^2=b^2=c=1$, $c\in {\bf c}$, and cyclic permutation. The number $m$ is called {\em the dihedral type} of $P$. 

Let $Q=\langle a, b, {\bf c} \ |\ (ab)^m,  a^2, b^2, c\in {\bf c} \rangle$. We have an epimorphism $\phi\colon G(P)\to G(Q)=D_m$. Let $\bar K(P)$ be the covering of $K(P)$ associated with the kernel.  Note that $\bar K(P)^{(1)}=\tilde K(Q)^{(1)}$, which is a $2m$-gon, consisting of double edges labeled in an alternating way with $a$ and $b$, and at every vertex we have a $c$ loop, for every $c\in {\bf c}$. An $a$-side of $\bar K(P)$ is a connected subgraph of the 1-skeleton that consists of a double edge with label $a$, together with all the $b$-double edges and $c$-loops connected to the two vertices of the $a$-double edge. A $b$-side is defined in an analogous way. We say $P$ is {\em side injective} if the inclusion of any side $S\to \bar K(P)$ is $\pi_1$-injective.  

\begin{lemma}\label{lem:long} Assume that $P$ is of dihedral type $m\ge 3$. Then $P$ is side injective if and only if every cyclically reduced word $w$ that represents the trivial element in $G(P)$ contains a reduced subword $s$ which is a cyclic permutation of 
$$a^{\alpha_1}d_1b^{\beta_1}d_2a^{\alpha_2}d_3b^{\beta_2}$$
or it's inverse.  The $\alpha_i$ and $\beta_i$ are odd integers and the $d_i$ are words in the generators containing $a$ and $b$ with even exponents (the $d_i$ could be trivial).
\end{lemma}

\noindent Proof. Assume first that $w$ is a cyclically reduced word that represents the trivial element in $G(P)$ and contains a reduced subword $s=a^{\alpha_1}d_1b^{\beta_1}d_2a^{\alpha_2}d_3b^{\beta_2}$ (the simplest setting is $s=abab$, a case the reader should have in mind). This subword does not lift into a side of $\bar K(P)$, and hence $w$ does not lift into a side. It follows that $P$ is side injective. 

Next assume that $P$ is not side injective. Assume $w$ is a cyclically reduced word that represents the trivial element of $G(P)$ and does lift into an $a$-side of $\bar K(P)$. Then $w$ is a reduced word in powers of $a^2$, $b^2$, $ab^2a^{-1}$, $c$, $aca^{-1}$, where $c\in {\bf c}$. So $w$ does not contain a subword of the form $s$. \qed 

\begin{example}\label{ex:torsion} $P=\langle a, b\ |\ (ab)^m\rangle$, $m\ge 3$, is side injective. This is because 1-relator presentations with torsion are Dehn presentations (in particular $G(P)$ is hyperbolic). See Newman \cite{Newman}. A word $w$ that is trivial in the group contains a subword of length more than $1/2$ of a cyclic permutation of the relator or its inverse, hence it contains a a cyclic permutation of $abab$, or its inverse. The result follows from Lemma \ref{lem:long}.
\end{example}

\begin{example} More generally, if $P=\langle a, b, {\bf c}\ |\ r(a,b, {\bf c}) \rangle$ (${\bf c}$ could be empty) is a Dehn presentation of dihedral type $m\ge 3$ so that more than half of a cyclic permutation of the relator or its inverse contains a subword $s$ as in Lemma \ref{lem:long}, then $P$ is side injective. Recall that $P$ is a Dehn presentation for instance in case it satisfies the small cancellation condition $C'(1/6)$ or $C'(1/4)-T(4)$ (see Chapter V, Theorem 4.4 in Lyndon and Schupp \cite{LS77}). For example if
$$r(a,b, {\bf c})=a^{\alpha_1}d_1b^{\beta_1}d_2a^{\alpha_2}d_3b^{\beta_2}d_4a^{\alpha_3}d_5b^{\beta_3}d_6a^{\alpha_4}d_7b^{\beta_4}d_8$$
where the $\alpha_i$ and $\beta_i$ are odd integers satisfying $|\alpha_i|=|\alpha_j|$, $|\beta_i|=|\beta_j|, \forall i,j\le 4$  and the $d_i$ are words of the same length containing $a$ and $b$ with even exponents, and $P$ satisfies the small cancellation condition $C'(1/6)$ or $C'(1/4)-T(4)$, then $P$ is side injective. Concrete examples are 
$$\langle a, b,c \ |\ (acbc^{-1}ac^{-1}bc)^2   \rangle$$
or
$$\langle a, b,c \ |\ acbc^{-1}acbcac^{-1}bc^{-1}ac^{-1}bc    \rangle$$
which are $C'(1/4)-T(4)$ and 
$$\langle a, b,c \ |\ acbca^{-1}cbc^{-1}a^{-1}c^{-1}bcac^{-1}bc^{-1}    \rangle$$
which is $C'(1/6)$. These presentations were checked with the help of {\sc GAP} (see \cite{GAP}) and the package {\sc SmallCancellation} by Ivan Sadofschi Costa (see \cite{ISC}).
\end{example}

\begin{example} The Artin presentation $P=\langle a, b\ |\ prod(a,b,m)=prod(b,a,m) \rangle$ is not side injective for $m=3$, but is side injective for $m\ge 4$:\\

\noindent 1) $m=3$. We show that $P=\langle a, b\ |\ aba=bab \rangle$ is not side injective. We have $a^2(aba^2ba)a^{-2}=aba^2ba$
in $G(P)$ because $(aba)^2=aba^2ba$ is central. So
$$w=a^2ba^2ba^{-2}b^{-1}a^{-2}b^{-1}=1$$
in $G(P)$. Note that $w$ lifts into a $b$-side of $\bar K(P)$.\\

\noindent 2) $m=4$. We show that $P=\langle a, b\ |\ abab=baba \rangle$ is side injective. 
Note that $x=abab$ is a central element. The quotient $G(P)/\langle x \rangle$ has a presentation $\langle a, b \ |\ (ab)^2 \rangle$. Let $y=ba$, then the presentation rewrites to $\langle a, y\ | y^2 \rangle$. In order to show that $P$ is $a$-side injective we have to show that $a^2$, $b^2$, $ab^2a^{-1}$ generate a free group of rank 3 in $G(P)$. We will do this by showing that $A=a^2$, $B=(ya^{-1})^2=ya^{-1}ya^{-1}$, and $C_0=a(ya^{-1})^2a^{-1}=aya^{-1}ya^{-1}a^{-1}$ generate a free group in the quotient presented by $Q=\langle a, y\ |\  y^2 \rangle=\mathbb Z*\mathbb Z_2$. Let $C_1=BC_0$. We have 
$$C_1=ya^{-1}ya^{-1}aya^{-1}ya^{-1}a^{-1}=ya^{-1}yya^{-1}ya^{-1}a^{-1}=ya^{-1}a^{-1}ya^{-1}a^{-1}=ya^{-2}ya^{-2}.$$ And finally let $C=C_1A=ya^{-2}y$. In summary we have 
$$A=a^2, \ B=ya^{-1}ya^{-1},\ C=ya^{-2}y.$$ The group $H=\langle A, B, C\rangle$ is a normal free subgroup of $G(Q)$ of rank 3 and index 4. Figure \ref{fig:4covering} shows a covering space $p\colon \bar K(Q)\to K(Q)$ so that $\pi_1(\bar K(Q))$ is free of rank 3 and $p_*(\pi_1(\bar K(Q)))=\langle A, B, C\rangle \le \pi_1(K(Q))$. The argument for $b$-side injectivity is analogous.
\begin{figure}[H] 
   \centering
   \includegraphics[width=2in]{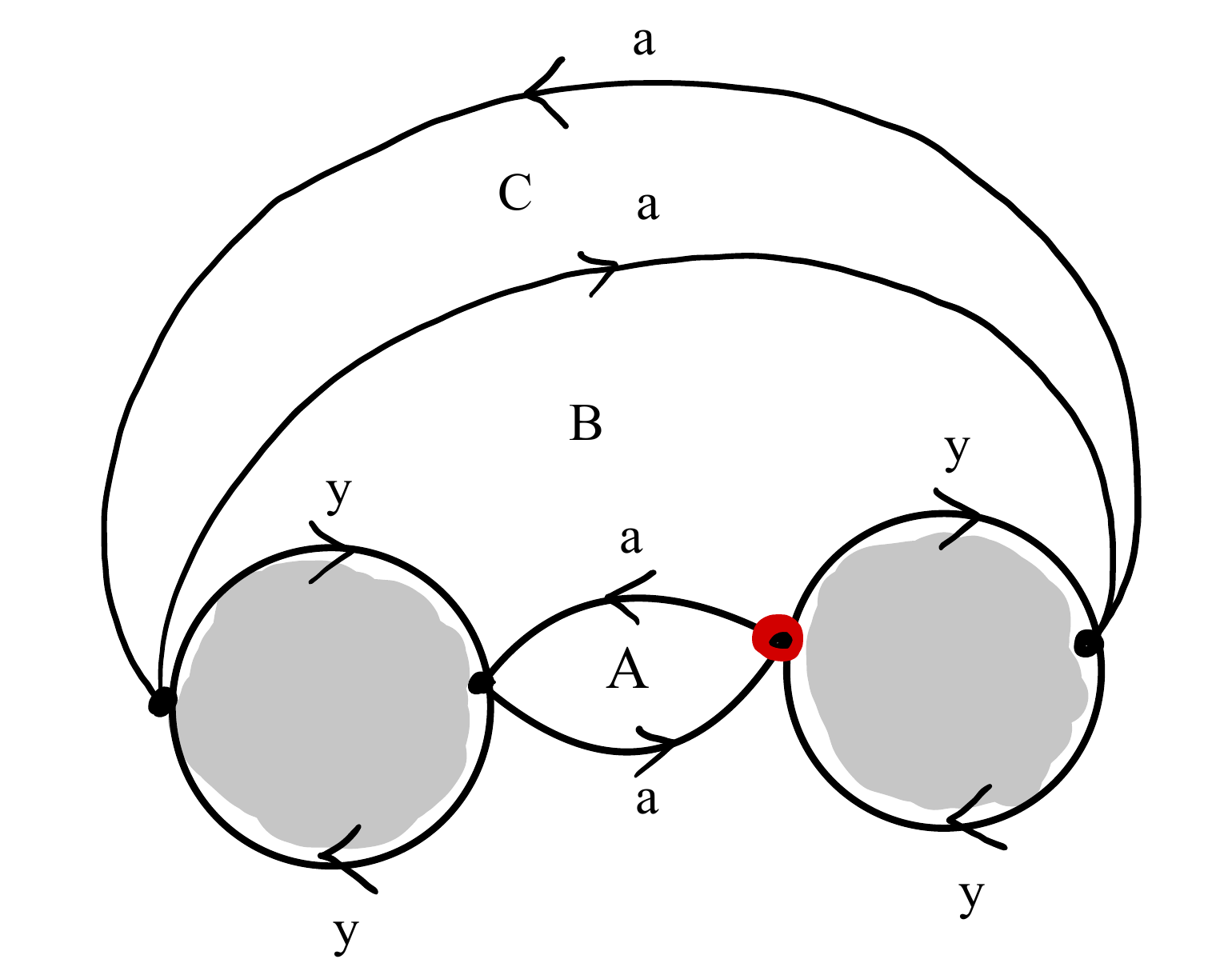} 
   \caption{If $Q=\langle a, y\ |\  y^2 \rangle$ then the universal covering $\tilde K(Q)$ is a tree with spheres attached. Here we see the intermediate covering $\bar K(Q)$ corresponding to the subgroup $H=\langle A, B, C\rangle$. The gray discs with boundary $y^2$ indicate 2-spheres.}
   \label{fig:4covering}
\end{figure}

\bigskip\noindent 3) $m\ge 6$ and even. This case is easy. Let $x=prod(a,b,m)$. The quotient $G(P)/\langle x \rangle$ is presented by $\langle a, b\ |\ (ab)^{\frac{m}{2}}\rangle$ which is a Dehn presentation, being a 1-relator presentation with torsion. Since $m\ge 6$, we have $\frac{m}{2}\ge 3$. Side injectivity follows from Example \ref{ex:torsion}.\\

\noindent 4) $m\ge 5$ and odd. Let $x=prod(a,b,m)$ and $y=ba$. Note that $x=ay^{\frac{m-1}{2}}$. Using $a=xy^{\frac{-m+1}{2}}$ and $b=y^{\frac{m+1}{2}}x^{-1}$ the presentation $P$ can be rewritten to $\langle x, y\ |\ x^2=y^m\rangle$. Thus $G(P)/\langle x^2\rangle$ is presented by $\langle x, y\ |\ x^2, y^m \rangle$ which is the hyperbolic group $\mathbb Z_2*\mathbb Z_m$. In the original generators this is $\langle a,b\ |\ prod(a,b,m)^2, (ba)^m\rangle$. If this were a Dehn presentation we could proceed as in the previous case (at least for $m\ge 7$), but we do not know. Instead we argue as in case 2. For simplicity we assume $m=5$, the other cases $m\ge 7$ go along the same lines. In order to show that $P$ is $a$-side injective we have to show that $a^2$, $b^2$, $ab^2a^{-1}$ generate a free group of rank 3 in $G(P)$. In terms of $x$ and $y$ it suffices to show that $xy^{-2}xy^{-2}$, $y^3x^{-1}y^3x^{-1}$, and 
$(xy^{-2})y^3x^{-1}y^3x^{-1}(xy^{-2})^{-1}$ generate a free subgroup of rank 3 in the quotient presented by $Q=\langle x, y\ |\ x^2, y^5 \rangle$. Let
\[ A=xy^3xy^3, \ B=y^3xy^3x, \ C_0=xyxy^3xy^2x.\] Let $C=C_0A=(xyxy^3xy^2x)(xy^3xy^3)=xyxy$.
$$A=xy^3xy^3, B=y^3xy^3x, C=xyxy.$$ Note that 
$$C(y^{-1}Cy)=(xyxy)y^{-1}(xyxy)y=xy^2xy^2=B^{-1}$$ and 
$$C(y^{-1}Cy)(y^{-2}Cy^2)=xy^2xy^2y^{-2}xyxyy^2=xy^3xy^3=A.$$ So it suffices that to show that 
$$X=C, \ Y=y^{-1}Cy, \ Z=y^{-2}Cy^2$$ generate a free subgroup of rank 3. 
Figure \ref{fig:5covering} shows a covering space $p\colon\bar K(Q)\to K(Q)$ so that $\pi_1(\bar K(Q))$ is free of rank 3 and $p_*(\pi_1(\bar K(Q)))=\langle X, Y, Z\rangle \le \pi_1(K(Q))$. The argument for $b$-side injectivity is analogous.
\begin{figure}[H] 
   \centering
   \includegraphics[width=3in]{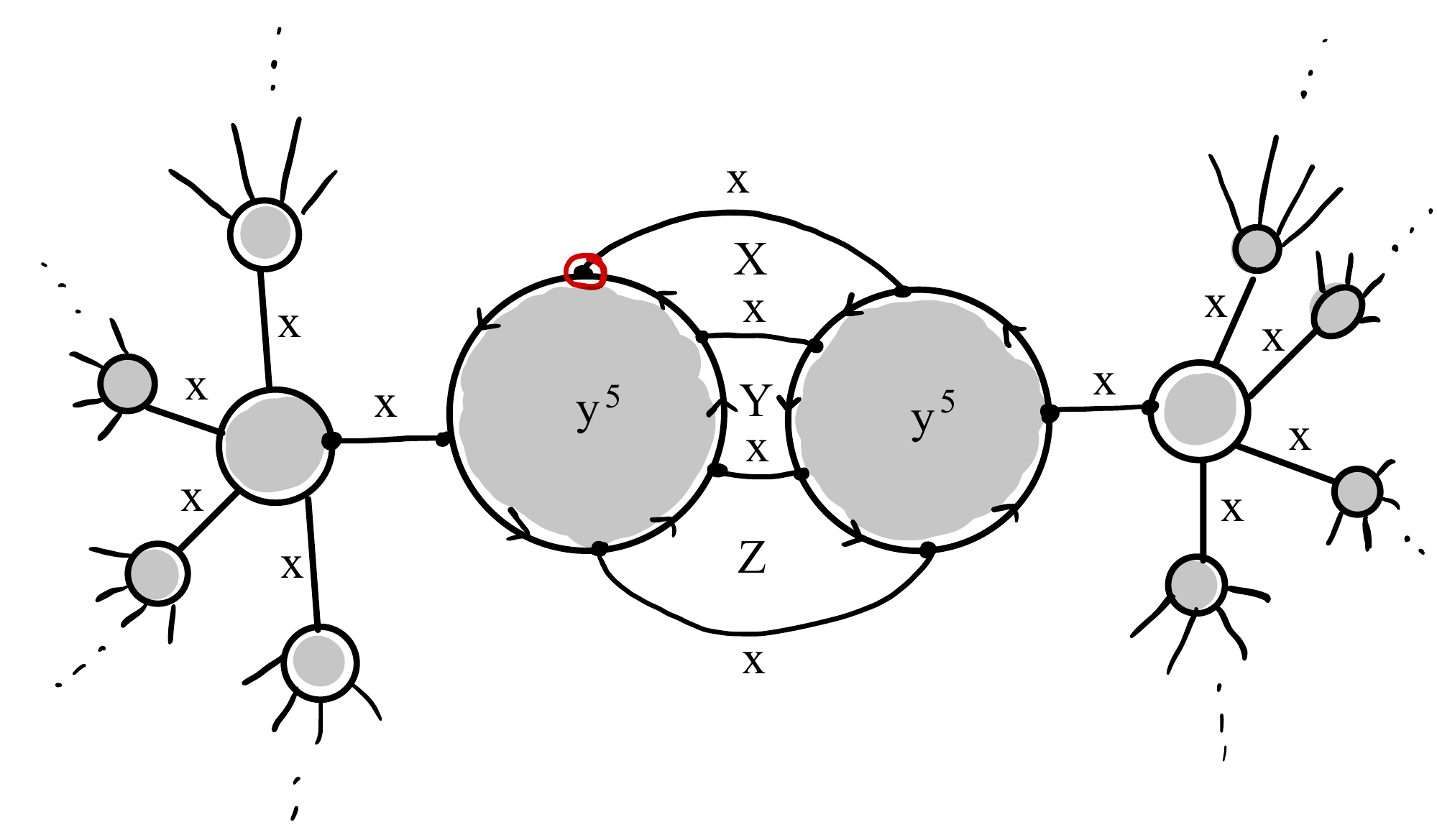} 
   \caption{A rendering of the covering space $\bar K(Q)$. Each $x$-edge represents a double $x$-edge into which 2 discs with boundary $x^2$ are glued. Each gray disc represents 5 discs with boundary $y^5$.}
   \label{fig:5covering}
\end{figure}
\end{example}

\begin{example} Let $P=\langle a, b, {c}\ |\ a(babcaba)=(babcaba)b \rangle$. Then $P$ is side injective by Theorem \ref{thm:ok}.
\end{example}

\begin{thm}\label{thm:ok} Suppose $P$ has dihedral type $m\ge 3$ and 
\begin{itemize}
\item $P=\langle a, b, {c}\ |\  a(u_1c^{\epsilon}u_3)=(u_1c^{\epsilon}u_3)b \rangle$, or
\item $P=\langle a, b, {\bf c}\ |\  a(u_1c_i^{\epsilon}u_2c_j^{\epsilon}u_3)=(u_1c_i^{\epsilon}u_2c_j^{\epsilon}u_3)b \rangle$, where
\end{itemize}

\begin{enumerate}
\item $c_i,c_j\in {\bf c}$ ($i=j$ is possible), $\epsilon=\pm 1$;
\item $u_1$ and $u_3$ do not contain any $c_k\in\bf c$ (or $c$ in the first case), and $u_2$ is arbitrary;
\item both $u_1^{-1}a$ and $u_3b^{-1}$ contain a subword $s$ as in Lemma \ref{lem:long}.
\end{enumerate}
Then $P$ is side injective.
\end{thm} 

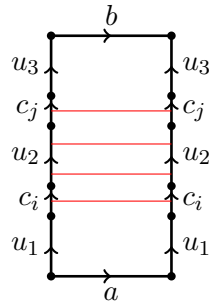
\begin{figure}[H]
   \centering
\begin{tikzpicture}[scale=0.8]
 
\fill (0,0) circle (2pt);
\fill (0,1) circle (2pt);
\fill (2,1) circle (2pt);
\fill (0,1.5) circle (2pt);
\fill (2,1.5) circle (2pt);
\fill (0,2.5) circle (2pt);
\fill (2,2.5) circle (2pt);
\fill (0,3) circle (2pt);
\fill (2,3) circle (2pt);
\fill (2,0) circle (2pt);
\fill (0,4) circle (2pt);
\fill (2,4) circle (2pt);
 
\begin{scope}[decoration={markings, mark=at position 0.5 with {\arrow{>}}}]
\draw [postaction={decorate},line width=1pt] (0,0) -- (2,0) node[midway, below]{$a$};
\draw [postaction={decorate},line width=1pt] (0,4) -- (2,4) node[midway, above]{$b$};
\draw [postaction={decorate},line width=1pt] (0,0) -- (0,1) node[midway, left]{$u_1$};
\draw [postaction={decorate},line width=1pt] (2,0) -- (2,1) node[midway, right]{$u_1$};
\draw [postaction={decorate},line width=1pt] (0,1.5) -- (0,2.5) node[midway, left]{$u_2$};
\draw [postaction={decorate},line width=1pt] (2,1.5) -- (2,2.5) node[midway, right]{$u_2$};
\draw [postaction={decorate},line width=1pt] (0,3) -- (0,4) node[midway, left]{$u_3$};
\draw [postaction={decorate},line width=1pt] (2,3) -- (2,4) node[midway, right]{$u_3$};
\end{scope}

\begin{scope}[decoration={markings, mark=at position 0.8 with {\arrow{>}}}]
\draw [postaction={decorate},line width=1pt] (0,1) -- (0,1.5) node[midway, left]{$c_i$};
\draw [postaction={decorate},line width=1pt] (2,1) -- (2,1.5) node[midway, right]{$c_i$};
\draw [postaction={decorate},line width=1pt] (0,2.5) -- (0,3) node[midway, left]{$c_j$};
\draw [postaction={decorate},line width=1pt] (2,2.5) -- (2,3) node[midway, right]{$c_j$};
\end{scope}
 
\draw[red] (0,1.25)--(2,1.25);
\draw[red] (0,1.7)--(2,1.7);
\draw[red] (0,2.2)--(2,2.2);
\draw[red] (0,2.75)--(2,2.75);

\end{tikzpicture}
   \caption{The relator disc drawn as a rectangle.}
   \label{fig:square}
\end{figure}

\begin{figure}[H]
   \centering
   \includegraphics[width=4.8cm, height=3.6cm]{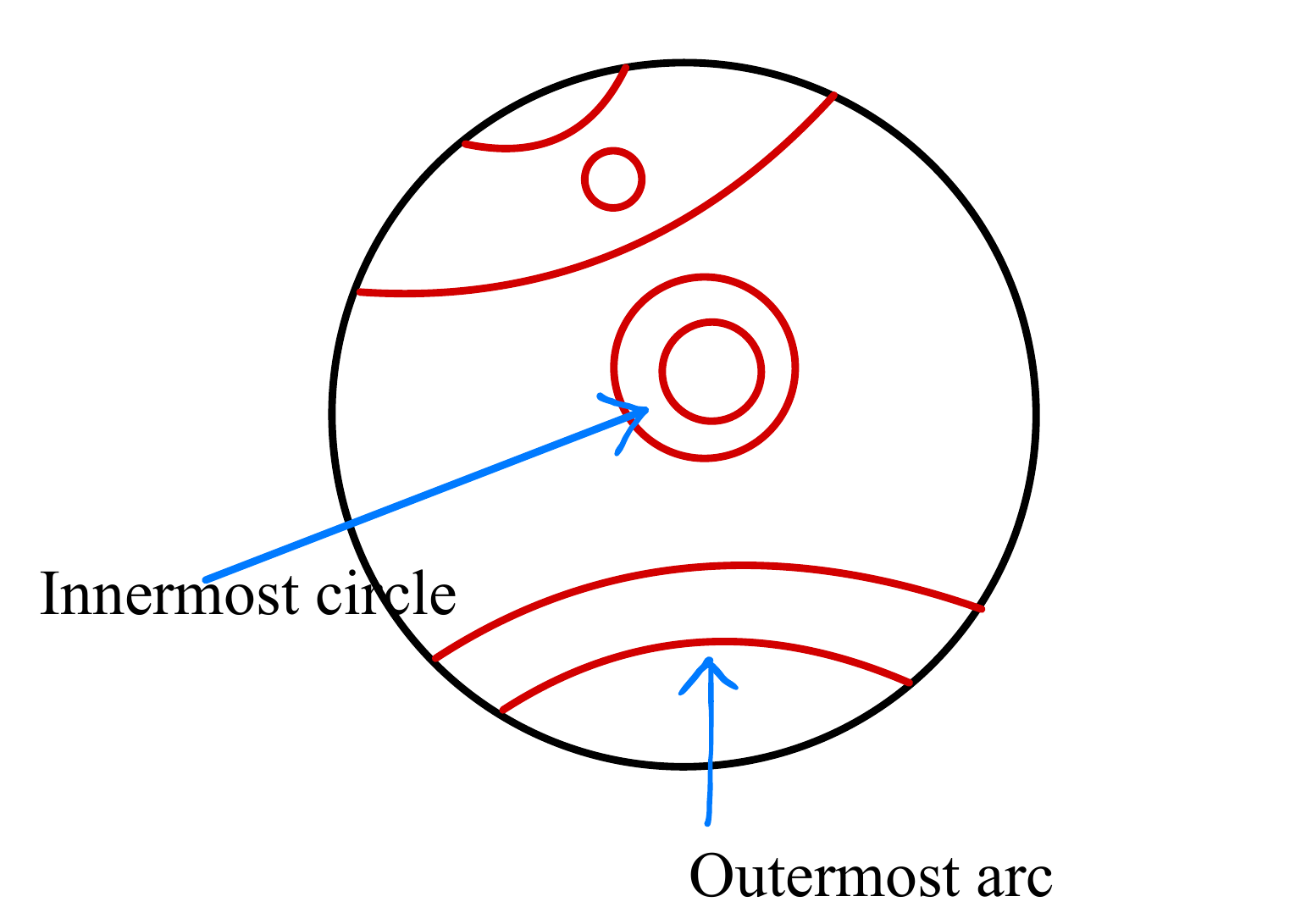} 
   \caption{A disc with red arcs, indicating innermost circles and outermost arcs.}
   \label{fig:redarcs}
\end{figure}

\noindent Proof. We assume we are in the second case and $\epsilon=1$. The first case is shown in an analogous way. Envision the relator disc placed in the plane as a rectangle, where the $a$ on the very left of the equation and the $b$ on the very right of the equation are horizontal edges, and the word $u_1c_iu_2c_ju_3$ is a vertical edge sequence. Connect the midpoints of $\bf c$-edges on the left and right by horizontal red edges. See Figure \ref{fig:square}.
Suppose that $w$ is a cyclically reduced word that represents the trivial element in $G(P)$. Let $D$ be a reduced Van Kampen diagram with boundary $w$. We may assume that $D$ is a topological disc. The red edges in our relator disc will form red circles and red arcs connecting points on the boundary of $D$. See Figure \ref{fig:redarcs}.
Consider an innermost red circle. Going around the inside we read off a word that freely reduces to $u_1^{-1}a^ku_1$ or $u_3b^ku_3^{-1}$, for some $k\in \mathbb Z$. If $k=0$, then $D$ is not reduced. If $k\ne 0$, then $G(P)$ has torsion. Both is not the case, hence there are no red circles in $D$. Consider an outermost red arc $\alpha$. Let $E$ be the component of $D-\alpha$ that does not contain anything red. Reading along the part of the boundary of $D$ which belongs to $E$ gives a reduced word (a subword of the reduced word $w$) equal to  $u_1^{-1}a^ku_1$ or $u_3b^ku_3^{-1}$. Because $D$ is reduced $k$ cannot be zero. If $k$ is positive then $u_1^{-1}a^ku_1$ contains $u_1^{-1}a$ and hence a word $s$ as in Lemma \ref{lem:long}. Also, $u_3b^ku_3^{-1}$ contains $bu_3^{-1}$, and since $u_3b^{-1}$ contains a word $s$ as in Lemma \ref{lem:long} so does 
$(u_3b^{-1})^{-1}=bu_3^{-1}$. The case where $k$ is negative goes the same way. It now follows from Lemma \ref{lem:long} that $P$ is side injective. \qed

\vspace{0.7cm}

\section{Last words about LOT applications}

\begin{thm} Let $\Gamma$ be a LOT of Coxeter type. Suppose that for every edge $e=(a\stackrel{w_e}{\to} b)$ the word $w_e$ is of the form $u_1c^{\epsilon}u_3$, or $u_1c^{\epsilon}u_2c^{\epsilon}u_3$ for some $c\ne a,b$, as in Theorem \ref{thm:ok}. Then $K(\Gamma)$ is aspherical.
\end{thm}

\noindent Proof. Each $\hat P_e$ is side injective. This follows from Theorem \ref{thm:ok}. Thus each $\bar K_e$ is side injective. The result follows from Theorem \ref{thm:aspherical}. \qed \\

\newpage

What if side injectivity fails?

\begin{thm} Suppose $\Gamma$ is a LOT of Coxeter type and there exist two edges $e_1$ and $e_2$ in $\Gamma$ so that 

\begin{enumerate}
\item $\bar K_{e_1}\cap \bar K_{e_2}=S$;
\item neither $\bar K_{e_1}$ or $\bar K_{e_2}$ is side injective, and in fact we have:\\
If $N_1=\ker (\pi_1(S)\to \pi_1(\bar K_{e_1}))$ and $N_2=\ker (\pi_1(S)\to \pi_1(\bar K_{e_2}))$ then $\frac{N_1\cap N_2}{[N_1,N_2]}\ne 1$.
\end{enumerate}
Then Whitehead's asphericity conjecture is false.
\end{thm}

\noindent Proof. Suppose Whitehead's conjecture is true. Then $K(\Gamma)$ and hence $\bar K(\Gamma)$ is aspherical. Note that $\bar K_{e_1}\cup \bar K_{e_2}$ is a subcomplex of $\bar K(\Gamma)$. Let $w$ be a reduced edge loop in $S$ that represents a non-trivial element in the quotient $\frac{N_1\cap N_2}{[N_1,N_2]}$. It is the boundary of a van Kampen diagram $D_1$ for $\bar K_{e_1}$ and also the boundary of a van Kampen diagram $D_2$ for $\bar K_{e_2}$.  The two diagrams can be glued together to form a non-trivial element in $\pi_2(\bar K_{e_1}\cup \bar K_{e_2})$ (see Gutierrez-Ratcliffe \cite{GutierrezRatcliffe}). A contradiction.\qed

\bigskip Jens Harlander, Boise State University

email: jensharlander@boisestate.edu

\bigskip Stephan Rosebrock, P\"adagogische Hochschule Karlsruhe, 

email: rosebrock@ph-karlsruhe.de


\begin{thebibliography}{99}
\bibitem{BerrickHillman} Berrick A.J., Hillman J.A.  {\em Whitehead's Asphericity Question and Its Relation to Other Open Problems}. In: Singh M., Song Y., Wu J. (eds) Algebraic Topology and Related Topics. Trends in Mathematics. Birkh\"auser, Singapore (2019).

\bibitem{Bleiler} S. A. Bleiler, {\em Two-generator cable knots are tunnel one}, Proceedings of the AMS, Volume 122, Number 4 (1994).

\bibitem{Bogley} W. A. Bogley, {\em J. H. C. Whitehead's asphericity question}. In 
Two-dimensional Homotopy and Combinatorial Group Theory (C. Hog-Angeloni, W. Metzler, and A. J. Sieradski, editors), Volume 197 of London Math. Soc. Lecture Note Ser. Cambridge University Press, 1993.

\bibitem{Button} J. O. Button, {\em Large groups of deficiency 1}, Isr. J. Math. 167, 111 (2008).

\bibitem{Davis} M. W. Davis, {\em The Geometry and Topology of Coxeter Groups}, LMS Monographs 32, Princeton University Press 2007.

\bibitem{CaretteWeidmann} M. Carette, R. Weidmann, {\em On the rank of Coxeter groups}, https://arxiv.org/pdf/0910.4997.pdf

\bibitem{GAP} GAP -- Groups, Algorithms, and Programming, Version 4.11.0, https://www.gap-system.org, (2020).

\bibitem{GutierrezRatcliffe} M. Gutierrez and J. Ratcliffe, On the second homotopy group, Quart. J. Math. Oxford (2) 32 (1981), 45-55.

\bibitem{HarRose} J. Harlander, S. Rosebrock, {\em Generalized knot complements and some
aspherical ribbon disc complements}, Journal of
Knot Theory and its Ramifications, Vol. 12, No. 7 (2003) 947-962.

\bibitem{HarRose2017} J. Harlander, S. Rosebrock, {\em Injective labeled oriented trees are aspherical}, Mathematische Zeitschrift 287 (1) (2017), p. 199-214.

\bibitem{HarRose2020}J. Harlander, S. Rosebrock, {\em Relative Vertex Asphericity}, to appear in the Bulletin of the Canadian Math. Soc. (2020), Math arXiv:1912.12512 [GT].

\bibitem{Howie} J. Howie, {\em On the asphericity of ribbon disc complements}, Transactions of the AMS, Vol. 289, No. 1 (1985) 281-302.

\bibitem{Hillman} J. A. Hillman, {\em Algebraic Invariants of Links}, 2nd edition, World Scientific Publ. Co. (2012).

\bibitem{KlimentoSakuma} E. Klimenko and M. Sakuma, {\em Two-generator discrete subgroups of Isom($\mathbb H^2$) containing orientation-reversing elements}, Geom. Dedicata 72 (1998), no. 3, 247-282.

\bibitem{LS77}
R.~Lyndon and P.~Schupp.
\newblock {\it Combinatorial group theory}, Springer Verlag, Berlin, (1977).

\bibitem{MulhollandRolfsen} J. Mulholland, D. Rolfsen, {\em Local indicability and commutator
subgroups of Artin groups}, https://arxiv.org/pdf/math/0606116.pdf

\bibitem{Newman} B. B. Newman, Some results on one-relator groups, Bull. Amer. Math. Soc. 74 (1968), 568-571.

\bibitem{Ro18} S. Rosebrock, {\em Labelled Oriented Trees and the Whitehead Conjecture}; 
Advances in Two-Dimensional Homotopy and Combinatorial Group Theory;
Cambridge University Press, LMS Lect.~Notes  446, editors
W.~Metzler, S.~Rosebrock; (2018), pp. 72--97.

\bibitem{ISC} I. Sadofschi Costa, {\em The small cancellation package}, https://github.com/isadofschi/smallcancellation

\bibitem{Serre} J. P. Serre, Trees, Springer (2003).

\bibitem{Weidmann} R. Weidmann, {\em The rank problem for sufficiently large Fuchsian groups}, Proc. LMS (3) 95 (2007), no. 3, pp. 609--652.

\bibitem{Yajima} T. Yajima, {\em On a characterization of knot groups of some knots in $\mathbb R^4$}, Osaka Math. J. 6 (1969), pp. 435-446.
\end{thebibliography}
\end{document}